\documentclass{article}
\usepackage{amsmath,amssymb,amscd,graphicx,psfrag,epsf,amsthm
}
\usepackage[all]{xy}

\newcommand{\PP}{\text {\bf P}}

\newcommand{\cC}{{\mathcal {C}}}
\newcommand{\cI}{{\mathcal {I}}}
\newcommand{\cT}{{\mathcal {T}}}

\newcommand{\cU}{{\mathcal {U}}}

\newcommand{\cE}{{\mathcal {E}}}

\newtheorem{thm}{Theorem}[section]

\newtheorem{lemma}[thm]{Lemma}
\newtheorem{cor}[thm]{Corollary}

\newtheorem{prop}[thm]{Proposition}

\newtheorem{conj}[thm]{Conjecture}
   
\theoremstyle{definition}
\newtheorem{defn}[thm]{Definition}

\newtheorem{notation}[thm]{Notation}   
  
\newtheorem{defn-thm}[thm]{Definition-Theorem}  

\theoremstyle{remark}

\renewcommand{\c}[0]{{\mathbb C}}  

\renewcommand{\o}[0]{{\mathcal O}}

\renewcommand{\a}[0]{{\mathbb A}}

\newcommand{\p}[0]{{\mathbb P}}

\newcommand{\map}[0]{\dasharrow}

\newcommand{\spec}[0]{\operatorname{Spec}}
\newcommand{\rat}[0]{\operatorname{RatCurves}^n}

\newcommand{\rank}[0]{\operatorname{rank}}

\newcommand{\codim}[0]{\operatorname{codim}}    
\newcommand{\im}[0]{\operatorname{im}}    
\newcommand{\proj}[0]{\operatorname{Proj}}    
    
\newcommand{\locus}[0]{\operatorname{locus}}

\newcommand{\Hom}[0]{\operatorname{Hom}}

\newcommand{\aut}[0]{\operatorname{Aut}}

\newcommand{\chow}[0]{\operatorname{Chow}}

\newcommand{\sym}[0]{\operatorname{Sym}}

\def\into{\DOTSB\lhook\joinrel\rightarrow}

\numberwithin{equation}{section}

\author{Carolina Araujo 
\\
  IMPA \\
  Rio de Janeiro, Brazil \\
  caraujo@impa.br
}
\title{Rational curves of minimal degree and characterizations of $\p^n$ }
\date{}

\begin{document}

\maketitle

\begin{abstract}
In this paper we investigate complex uniruled varieties $X$ whose rational curves of minimal
degree satisfy a special property.
Namely, we assume that the tangent directions to such curves at a general point $x\in X$ 
form a linear subspace of $T_xX$.
As an application of our main result, we give a unified geometric proof of Mori's,
Wahl's, Campana-Peternell's
and Andreatta-Wi\'sniewski's characterizations of $\p^n$.
\end{abstract}


\section{Introduction}\label{introduction}

Let $X$ be a smooth complex projective variety, and assume that $X$
is uniruled, i.e., there exists a rational curve through every point of $X$.
Let $H$ be a covering family of rational curves on $X$ having minimal 
degree with respect to some fixed ample line bundle. 
For each $x\in X$ denote by $\cC_x$ the subvariety of the 
projectivized tangent space at $x$ consisting of 
tangent directions to rational curves from $H$ passing through $x$. 
We are interested in varieties $X$ for which $\cC_x$ is a linear subspace
of $\p(T_xX)$ for general $x\in X$.
We prove the following result.

\begin{thm}\label{main_intro}
Suppose $\cC_x$ is a $d$-dimensional linear subspace of $\p(T_xX)$ for a general point $x\in X$.
Then there is a dense open subset $X^0$ of $X$ and a $\p^{d+1}$-bundle 
$\varphi^0:X^0\to T^0$ such that any curve from $H$ meeting $X^0$
is a line on a fiber of $\varphi^0$.
\end{thm}

In fact we prove a stronger result, as we allow $\cC_x$ to be a union 
of linear subspaces of $\p(T_xX)$ for general $x\in X$ 
(see Theorem~\ref{Cx_fibration}).
We remark that the variety $\cC_x$ 
has been studied in a series of papers by Hwang and Mok (see \cite{hwang}).

As an application, we provide a unified geometric proof of
the following characterization of $\p^n$.

\begin{thm}\label{p^n}
Let $X$ be a smooth complex projective $n$-dimensional variety.
Assume that the tangent bundle $T_X$ contains 
an ample locally free subsheaf $E$ of rank $r$. 
Then $X\cong \p^n$ and either $E\cong \o_{\p^n}(1)^{\oplus r}$ or
$r=n$ and $E= T_{\p^n}$.
\end{thm}

The first instance of this theorem, namely the case $E\cong T_X$,
was proved by Mori in \cite{mori79}. In his proof, Mori recovers the projective space
by studying rational curves of minimal degree passing through a general point of $X$. 
Then, in \cite{wahl}, Wahl settled the case that $E$ is a line bundle.
Wahl's proof is very different from Mori's. 
It relies on the theory of algebraic derivations in characteristic zero.
It does not make any use of the geometry of rational curves on $X$.
Recently Druel gave a geometric proof of Wahl's theorem in \cite{druel}.
His proof is based on studying the foliation by curves defined by the inclusion 
$E\into T_X$, and applying a criterion for algebraicity of the leaves.
In \cite{campana_peternell}, Campana and Peternell proved the theorem in the cases
$r=n$, $n-1$ and $n-2$.
The proof was finally completed by  Andreatta and Wi\'sniewski in \cite{andreatta_wisniewski}.
Their proof uses the geometry of minimal covering family of rational 
curves on $X$. It relies, on one side, on Mori's theorem, 
and on the other side, on Wahl's theorem.

Our proof follows the lines of Mori's proof of the Hartshorne conjecture in \cite{mori79}.
Here is the outline.
An $n$-dimensional variety $X$ whose 
tangent bundle contains an ample locally
free subsheaf is uniruled.
So we fix a covering family of rational curves of minimal degree on $X$, 
and consider the variety $\cC_x\subset \p(T_xX)$ of tangent directions to 
curves passing through $x\in X$.
We translate the existence of an ample locally free 
subsheaf of $T_X$ 
into projective properties of the embedding $\cC_x\into \p(T_xX)$,
and show that $\cC_x$ is a linear subspace of $\p(T_xX)$ for general $x\in X$.
By Theorem~\ref{main_intro}, there is a $\p^{d+1}$-fibration on 
a dense open subset of $X$, which can be extended in codimension $1$ following an argument in \cite{andreatta_wisniewski}.
If $d+1<n$, then the relative tangent bundle of 
such fibration inherits the ampleness
properties of $T_X$. 
We reach a contradiction by applying a result by Campana and Peternell
(\cite{campana_peternell}).

Throghout the paper we work over $\c$. 
In section~\ref{section:C_x} we gather some properties of minimal families of rational 
curves and the embedding $\cC_x\into \p(T_xX)$.
In section~\ref{section:linear_Cx} we investigate  varieties $X$ for which $\cC_x$ is 
a union of linear subspaces of $\p(T_xX)$ for general $x\in X$.
In section~\ref{section:proof} we give a unified proof of Theorem~\ref{p^n}.

\medskip

\noindent {\it Notation. } In our discussion on rational curves we follow the notation of \cite{kollar}.
By a general point of a variety $X$, we mean a point in some dense open subset of $X$.
If $E$ is a vector bundle on a variety $X$, we denote by $\PP(E)$
the Grothendieck projectivization $\proj_X(\sym (E))$.
If $V$ is a complex vector space,
we denote by $\p(V)$ the natural projectivization of $V$.
(So $\PP(V)=\p(V^{\vee})$.)


\section{Tangent directions to rational curves of minimal degree}\label{section:C_x}

Let $X$ be a smooth complex projective variety, and assume that $X$
is uniruled.
Let $H$ be an irreducible component of $\rat(X)$. 
We say that $H$ is a \emph{covering family} if the 
corresponding universal family dominates $X$. 
A covering family $H$ of rational curves on $X$ is called
\emph{minimal} if, for a 
general point $x\in X$, the subfamily of $H$ parametrizing curves through 
$x$ is proper.
It is called \emph{unsplit} if $H$ itself is proper.

Fix a minimal covering family $H$
of rational curves on $X$ (for instance, one can take $H$ to be a 
covering family having minimal degree with respect to some fixed ample
line bundle on $X$).

Let $x\in X$ be a general point and denote by $H_x$ the normalization of the 
subscheme of $H$ parametrizing rational curves passing through $x$.
Let $\pi_x:U_x\to H_x$ and $\eta_x:U_x\to X$
be the universal family morphisms,
\[
\begin{CD}
  U_x @>\eta_x >> X, \\
  @V\pi_x VV \\
  H_x
\end{CD}
\]
so that $U_x$ is normal and $\pi_x$ is a $\p^1$-bundle (see \cite[II.2.12]{kollar}). 
Denote by $\locus(H_x)$ the closure of the image of $\eta_x$
(with the reduced scheme structure).
We remark that a rational curve smooth at $x$ is parametrized by at most one element of $H_x$.


\begin{notation}
Let $f:\p^1\to X$ be a morphism birational onto its image  
such that $f(o)=x$. We denote by 
$[f]$ the element of $V$ (or $V_x$) parametrizing $f$.
Sometimes we also denote by $[f]$ the point $\varphi_x([f])\in H_x$ parametrizing the image of $f$.
It should be clear from the context whether we view $[f]$ as a member
of $\Hom(\p^1,X)$ or $\rat(X)$.
\end{notation}

Next we gather some important properties of minimal covering families of rational 
curves.

\begin{prop}\label{properties}
Let the notation be as above and $x\in X$ a general point . 
\begin{enumerate}
   \item  For every $[f]\in H_x$,
     $f^*T_X\cong \bigoplus\limits_{i=1}^n \o(a_i)$, 
     with all $a_i\geq 0$.
   \item  $H_x$ is a smooth projective variety of dimension $d:=\deg(f^*T_X)-2$. 
   \item  If $[f]$ is a general member of any irreducible component of $H_x$, then
     $
     f^*T_X\cong \o(2)\oplus \o(1)^{\oplus d}\oplus \o^{n-d-1}. 
     $
   \item  Every curve parametrized by $H_x$ is immersed at $x$ 
     (i.e., $df_o$ is nonzero for every $o\in \p^1$ such that $f(o)=x$). 
   \item  The dimension of the subscheme of $H_x$ parametrizing curves singular at $x$ 
     is at most the dimension of the subscheme of $H_x$ parametrizing curves with 
     cuspidal singularities (not necessarily at $x$).
   \item  If all the curves parametrized by $H_x$ are smooth at $x$,
     then the restriction of $\eta_x$ to each irreducible component of $U_x$
     is birational onto its image.
\end{enumerate}
\end{prop}

\begin{proof}
Property (1) follows from \cite[II.3.11]{kollar} and the assumption that $x$ is a general point.
Property (2) follows from \cite[II.1.7, II.2.16]{kollar} and the assumption that $H$ is a minimal covering family.
Property (3) follows from \cite[IV.2.9]{kollar}.

Properties (4) and (5) can be found in \cite{kebekus}. 
Property (5) is not explicitly stated in \cite{kebekus}, but follows
from the proof of \cite[Theorem 3.3]{kebekus}.

Property (6) is due to Miyaoka (see \cite[V.3.7.5]{kollar}).
\end{proof}

\begin{defn}
Define the tangent map 
$
\tau_x:  \ H_x  \map \p(T_xX)
$
by sending a curve that is smooth at $x$ to its tangent 
direction at $x$.

Define $\cC_x$ to be the closure of the image of $\tau_x$
in $\p(T_xX)$.
\end{defn}

\begin{thm} \label{tau}
Let the notation be as above. Then
\begin{enumerate}
\item {\bf (\cite{kebekus})} $\tau_x:H_x\to \cC_x$ is a finite morphism,
\item {\bf (\cite{hwang_mok_birationality})} $\tau_x:H_x\to \cC_x$ is birational, and thus
\item $\tau_x:H_x\to \cC_x$ is the normalization.
\end{enumerate}
\end{thm}

Notice that $\cC_x$ comes with a natural projective 
embedding into $\p(T_xX)$. 
It turns out that, for a general member $[f]\in H_x$,
the ``positive'' directions of $f^*T_X$ at $x$ determine
the tangent space of $\cC_x$ at $\tau_x([f])$. 
This is made precise in the next proposition.

\begin{defn}
Let $f:\p^1 \to X$ be a morphism birational onto its image such that $x=f(o)$.
Define the \emph{positive tangent space at $x\in X$ with respect to $f$} 
to be the following linear subspace of $T_xX$:
$$
T_xX^+_f \ := \ \im[H^0(\p^1,f^*T_X(-1))\to (f^*T_X(-1))_o\cong T_xX].  
$$
\end{defn}

\begin{prop} \label{TC_x}
Let $[f]\in H_x$ be a general element.
Then $\p(T_xX^+_f)\subset \p(T_xX)$ is the projective tangent space of 
$\cC_x$ at $\tau_x([f])$.
\end{prop}

\begin{proof}
See \cite[Proposition 2.3]{hwang} or \cite[Lemma 2.1]{andreatta_wisniewski}.
\end{proof}

From the splitting type of $f^*T_X$ one can check whether
$\tau_x$ is an immersion at $[f]$.

\begin{prop} \label{immersion}
The morphism $\tau_x$ is an immersion at $[f]\in H_x$ if and only if 
$
f^*T_X\cong \o_{\p^1}(2)\oplus \o_{\p^1}(1)^{\oplus d} \oplus  \o_{\p^1}^{\oplus n-d-1},
$
where $d=\deg (f^*T_X)-2$.
\end{prop}

\begin{proof}
Let $H'_x$ be an irreducible component of $H_x$. Let $V_x$ be the corresponding irreducible 
component of $\Hom(\p^1,X,o\mapsto x)$,
i.e., $V_x$ parametrizes morphisms (birational onto their images) 
whose images are parametrized by $H'_x$.
By Proposition~\ref{properties}(4),
every morphism parametrized by $V_x$ is an immersion at $o$. So we can define the morphism
$\cT_x: V_x  \to \p(T_xX)$ by setting $\cT_x([f])= \p(df_o(T_o\p^1))$.

We have the following commutative diagram:
\[
\xymatrix{
V_x \ar[d]_{\varphi_x} \ar[r]^{\cT_x}  & \ \p(T_xX), \\
H'_x \ar[ru]_{\tau_x} 
}
\]
where $\varphi_x$ is a smooth morphism with
fibers isomorphic to $\aut(\p^1,o)$.

Fix $[f]\in V_x$ and write 
$f^*T_X\cong \bigoplus\limits_{i=1}^n \o(a_i)$, 
with $a_1\geq \dots \geq a_n\geq 0$, $a_1\geq 2$, and
$\sum\limits_{i=1}^{n}a_i=2+d$.
Let us describe the tangent map 
$d\cT_x([f]):T_{[f]}V_x\to T_{\cT_x([f])}\p(T_xX)$ explicitly.

There are isomorphisms 
$T_{[f]}V_x \ \cong \ H^0(\p^1,f^*T_X\otimes \cI_o)$,
where $\cI_o$ denotes the ideal sheaf of $o$ in $\p^1$
(see \cite[II.1.7]{kollar}),  
and $T_{\cT_x([f])}\p(T_xX) \ \cong \ T_xX/\hat \cT_x([f])$,
where  $\hat \cT_x([f])$ denotes 
the $1$-dimensional subspace of $T_xX$ corresponding to the point 
$\cT_x([f])\in \p(T_xX)$.

Fix a local parameter $t$ of the local ring of $o$ on $\p^1$. 
If $v$ is a global section of $f^*T_X$ vanishing at $o$, then $d\cT_x([f])(v)$ 
is given by
$$
d\cT_x([f])(v)=\left[ \left. \frac{d}{dt}v(t)\right|_{t=o}\right]\in 
(f^*T_X)_o/T_o\p^1 \cong T_xX/\hat \cT_x([f]).
$$

So we see that $\im d\cT_x([f])\cong T_xX^+_f/\hat \cT_x([f])$.

Set $l= \dim \im d\cT_x([f]) = \sharp \{a_i|a_i>0\}-1$. 
Then $l\leq d=\sum\limits_{i=1}^{n}a_i - 2$, and 
equality holds if and only if 
$\tau_x$ is an immersion at $\varphi_x([f])$.
Since $a_1\geq 2$ and 
$a_i \geq 1$ for $1\leq i\leq l+1$, this is equivalent to 
$a_1=2$ and $a_i=1$ for $1<i\leq d+1$,
i.e., $f^*T_X\cong \o(2)\oplus \o(1)^{\oplus d} \oplus \o^{\oplus n-d-1}$.
\end{proof}

\begin{cor} \label{Cx_smooth}
If every irreducible component of $\cC_x$ is smooth, then  
\begin{enumerate}
   \item all curves parametrized by $H_x$ are smooth at $x$, and
   \item the restriction of the universal family morphism $\eta_x:U_x\to X$ 
         to each irreducible component of $U_x$ is birational onto its image.
\end{enumerate}
\end{cor}

\begin{proof}
Since every irreducible component of $\cC_x$ is smooth, 
$\tau_x$ is an immersion by Theorem~\ref{tau} 
(in fact, the restriction of $\tau_x$ to each irreducible component 
of $H_x$ is an isomorphism).
Thus, by Proposition~\ref{immersion},  
$f^*T_X=\o_{\p^1}(2)\oplus \o_{\p^1}(1)^{\oplus d} \oplus  \o_{\p^1}^{\oplus n-d-1}$
for every member $[f]\in H_x$.
From the splitting type of $f^*T_X$ we see that 
no curve parametrized by $H_x$ has a cuspidal singularity.
The corollary then follows from Proposition~\ref{properties}(5)--(6).
\end{proof}


\section{The distribution defined by linear $\cC_x$}\label{section:linear_Cx}

Let $X$ be a smooth uniruled complex projective variety. 
Let $H$ be a  minimal covering family 
of rational curves on $X$, and let $\cC_x$ be the subvariety of $\p(T_xX)$ defined in 
section~\ref{section:C_x}.
In this section we study varieties $X$ for which $\cC_x$ is 
a union of linear subspaces of $\p(T_xX)$ for general $x\in X$.

Denote by $\pi:  U\to H$ and $\eta: U \to X$ 
the universal family morphisms.
Consider the Stein factorization of $\eta$:
\[
\xymatrix{
                                      &  &    X' \ar[d]^{\rho} \\
U \ar[rr]^{\eta} \ar@/^0.3cm/[rru]^{\eta'} \ar[d]_{\pi}   & &     X. \\
H
}
\]
We may view $H$ as a minimal
covering family of rational curves on $X'$.

The main result in this section is the following.

\begin{thm}\label{Cx_fibration}
Suppose that $\cC_x$ is 
a union of $d$-dimensional linear subspaces of $\p(T_xX)$ for general $x\in X$.
Let $X'$ be as defined above. 
Then there is a dense open subset $U^0$ of $X'$ and a $\p^{d+1}$-bundle 
$\varphi^0:U^0\to T^0$ such that any curve on $X'$ parametrized by 
$H$ and meeting $U^0$
is a line on a fiber of $\varphi^0$.
\end{thm}

For a general point $x\in X$,
denote by $H_x^i$, $1\leq i\leq k$, the irreducible components of $H_x$,
and by $\cC_x^i$ the image of $H_x^i$ under $\tau_x$.
Suppose that each $\cC_x^i$ is a $d$-dimensional linear subspace of $\p(T_xX)$.

Viewing $H$ as a minimal covering family of rational curves on $X'$,
$H_{x'}$ is irreducible and $\cC_{x'}$ is a linear subspace of $\p(T_{x'}X')$
for a general point $x'\in X'$.
Moreover $X'$ is smooth along $\locus (H_{x'})$ 
(for $\eta$ is smooth along $\pi^{-1}(H_{x'})$ by \cite[II.3.5.3, II.2.15]{kollar}).

We obtain a rank $d+1$ distribution $D$ on a dense open subset of $X'$ as follows. 
For a general point $x'\in X'$, set
$D_{x'}=\hat \cC_{x'}$, the linear subspace of $T_{x'}X'$ corresponding 
to $\cC_{x'}\subset \p(T_{x'}X')$.

\begin{lemma}\label{Cx_integrable}
Suppose that $\cC_x$ is 
a union of linear subspaces of $\p(T_xX)$ for general $x\in X$.
Let $X'$ and $D$ be as defined above. 
Then the distribution $D$ is tangent to $\locus(H_{x'})$ for a general point $x'\in X'$.
In particular, $\locus (H_{x'})$ is smooth at $x'$.
\end{lemma}

\begin{proof}
Let $x'\in X'$ be a general point and set
$Y:= \locus(H_{x'})$.
By Frobenius' Theorem,
we are done if we can show that $T_{y}Y= D_{y} = \hat \cC_{y}$
for a general point $y\in Y$.

Let $[f]\in H_{x'}$ be a general member and let $y$ be a general point
in the image of $f$.
Let $o,p\in \p^1$ be such that $f(o)=x'$ and  $f(p)=y$.
Let $V_{x'}$ be the irreducible component of $\Hom(\p^1,X,o\mapsto x')$
corresponding to $H_{x'}$.
We have the following commutative diagram:

\[
\xymatrix{
\p^1 \times V_{x'} \ar[d] \ar[r]\ar@/^0.7cm/[rr]^{F} & U_{x'} \ar[d]^{\pi_{x'}} \ar[r]^{\eta_{x'}} & X'. \\
V_{x'} \ar[r]^{\varphi_{x'}} & H_{x'}
} 
\]

By generic smoothness,
the tangent space $T_{y}Y$ is the image in $T_yX'$
of the differential $dF_{(p,[f])}$.
From the description of
$dF_{(p,[f])}$ given in \cite[II.3.4]{kollar}, together with Proposition~\ref{TC_x},
we see that 
this is precisely $T_y{X'}^+_f = \hat \cC_{y}$. 
\end{proof}

Next we describe $\locus (H_{x'})$ for general $x'\in X'$.

\begin{lemma} \label{leaf_p^n}
Suppose that $\cC_x$ is 
a union of $d$-dimensional linear subspaces of $\p(T_xX)$ for general $x\in X$.
Let $X'$ be as defined above. 
Then, for a general point $x'\in X'$, 
the normalization of $\locus (H_{x'})$ is isomorphic to $\p^{d+1}$.
Under this isomorphism, the rational curves on $\locus (H_{x'})$ parametrized by $H_{x'}$
come from lines on $\p^{d+1}$ passing through a fixed point.
\end{lemma}

\begin{proof}
Set $Y=\locus(H_{x'})$ and let
$n:\tilde Y\to Y$ be the normalization.

The subfamily $H_Y=\{[f]\in H|f(\p^1)\subset Y\}$ is a minimal
covering family of rational curves on $Y$. 
Moreover, $x'$ is a general point of $Y$ (indeed,
$Y= \locus(H_{y})$ for a general point $y\in Y$ by Lemma~\ref{Cx_integrable}).
The subfamily $H_{Y,x'}$ of $H_Y$ 
parametrizing curves through $x'$ is just $H_{x'}\cong \p^d$.

Denote by $\pi_{x'}:  U_{x'}\to H_{Y,x'}$ and $\eta_{x'}: U_{x'} \to Y$ 
the universal family morphisms.
Since $\pi_{x'}$ is a $\p^1$-bundle, $U_{x'}$ is smooth.
Moreover, $\eta_{x'}$ is birational by Proposition~\ref{properties}(6).
We have the commutative diagram
\[
\xymatrix{
                                      &  &    \tilde Y \ar[d]^{n} \\
U_{x'} \ar[rr]^{\eta_{x'}} \ar@/^0.3cm/[rru]^{\tilde \eta_{x'}} \ar[d]_{\pi_{x'}}   & &     Y. \\
H_{Y,x'}
}
\]

Since $Y$ is smooth at $x'$, there is a unique point $\tilde x\in \tilde Y$
such that $n(\tilde x)=x'$, and $\tilde Y$ is smooth at $\tilde x$.

Let $\sigma \subset U_{x'}$ be the section of $\pi_{x'}$ that is contracted 
to $\tilde x \in \tilde Y$ by $\tilde \eta_{x'}$. 
Then $\tilde \eta_{x'}:U_{x'}\to \tilde Y$ is a surjective birational
morphism and restricts to an isomorphism on $U_{x'}\setminus \sigma$.
In particular $\tilde Y$ is smooth. 
In this setting, a standard argument by Mori (see \cite[V.3.7.8]{kollar}) yields
the result.
\end{proof}

\begin{proof}[Proof of Theorem~\ref{Cx_fibration}]
Let $X'$ be as defined above. 
By Lemmas \ref{Cx_integrable} and \ref{leaf_p^n}, 
together with Frobenius' Theorem,
there exists a dense open subset $U^0\subset X'$ and a morphism 
$\varphi^0:U^0\to T^0$ such that the normalization of the closure 
of the general fiber of $\varphi^0$ is isomorphic to $\p^{d+1}$.
By enlarging $U^0$ if necessary, we may assume that $X'\setminus U^0$
is the indeterminacy locus of $\varphi^0$.
If $d=n-1$, there is nothing to prove. So we assume that $\dim T^0\geq 1$.

Suppose that the general fiber of $\varphi^0$ is not proper (so that its closure
intersects $X'\setminus U^0$). 
Let $t\in T^0$ be a general point.
Then there exists a point $y\in X'\setminus U^0$ and a 
positive dimensional irreducible subvariety $T'\subset T^0$ containing $t$
such that $y$ lies in the closure of every fiber of $\varphi^0$ over $T'$.
Let  $H_y$ be the subscheme of $H$ parametrizing curves passing through $y$.

Let $t'\in T'$ be a general point, and $x'$ a general point in the fiber over $t'$.
By Lemma~\ref{leaf_p^n}, the normalization of $\locus (H_{x'})$ is isomorphic to $\p^{d+1}$, 
and the curves parametrized by $H_{x'}$ come from lines in $\p^{d+1}$.
This has two consequences. 
First, there is an element $[f]\in H_{x'}$ parametrizing a curve passing through $y$.
Since $x'$ is general, 
$f^*T_{X'}\cong \o_{\p^1}(2)\oplus \o_{\p^1}(1)^{\oplus d} \oplus  \o_{\p^1}^{\oplus n-d-1}$,
and thus $\dim_{[f]}H_y=d$ by \cite[II.1.7, II.2.16]{kollar}.
Second, $\locus (H_{x'})\subset \locus (H_y)$.
Since this holds for a point $x'$ in a general fiber over $T'$, we have that 
$\dim_{[f]} H_y\geq d+1$, contradicting the equality obtained above.

We conclude that the general fiber of $\varphi^0$ is proper. By shrinking $U^0$ and $T^0$ 
if necessary we get that $\varphi^0:U^0\to T^0$ is a $\p^{d+1}$-bundle.
\end{proof}

When $H_x$ is irreducible and $\cC_x$ is a linear subspace of $\p(T_xX)$,
Theorem~\ref{Cx_fibration} yields a dense open subset $X^0\subset X$ and a 
$\p^{d+1}$-bundle $\varphi^0:X^0\to T^0$. 
If we further assume that $H$ is an unsplit family, then $\varphi^0$
can be extended in codimension $1$, as we show below.

\begin{thm}\label{extending_in_codim_1}
Suppose $H$ is an unsplit family and 
$\cC_x$ is a linear subspace of $\p(T_xX)$ for general $x\in X$.
Then there is an open subset $X^0\subset X$ whose complement has codimension at least $2$ 
in $X$, and a $\p^{d+1}$-bundle 
$\varphi^0:X^0\to T^0$ over a smooth base 
satisfying the following property.
Every rational curve parametrized by $H$ and meeting $X^0$
is a line on a fiber of $\varphi^0$.
\end{thm}

\begin{proof}
We follow an argument in \cite{andreatta_wisniewski}.

Let $\varphi^0:X^0\to T^0$ be the $\p^{d+1}$-bundle from Theorem~\ref{Cx_fibration}.
Let $T\to \chow(X)$ be the normalization of the 
closure of the image of $T^0$ in $\chow(X)$, 
and let $\cU$ be the normalization of the universal family 
over $T$. Denote by $p:\cU\to T$ and $q:\cU\to X$ 
the universal family morphisms.

Let $0\in T$ be any point.
Set $\cU_0=p^{-1}(0)$ and 
let $x,y\in \cU_0$ be arbitrary points.
We can find a $1$-parameter family of fibers $\cU_t=p^{-1}(t)$,
together with points $x_t,y_t\in \cU_t$,
such that $\cU_t\cong \p^{d+1}$ for $t\neq 0$,
and $\lim_{t\to 0}(\cU_t,x_t,y_t)=(\cU_0,x,y)$.

Let $l_t\subset \cU_t$ be the curve parametrized by $H$ 
joining $x_t$ and $y_t$. 
Since $H$ is unsplit, the limit $\lim_{t\to 0}[l_t]$ lies
in $H$.
It parametrizes an irreduclible (and reduced) rational curve 
$l\subset \cU_0$ joining $x$ and $y$.
This shows that $\cU_0$ is irreducible.

Notice that $q:\cU\to X$ is birational and $T$ has dimension $n-d-1$.
Let $E\subset \cU$ be an irreducible component of the exceptional 
locus $E'$ of $q$. Since $X$ is smooth, $E$ has codimension $1$ in $\cU$.
Set $p_E=p|_E$ and $\cE=p(E)\subset T$.

Let $\cU_t$ be an arbitrary fiber of $p$ and assume that 
$\cU_t\cap E\neq \varnothing$, i.e., $t\in \cE$.
Since $E$ misses the general fiber of $p$, 
$\dim \cE \leq \dim T -1=n-d-2$.
Set $E_t=p_E^{-1}(t)$. Then 
$d+1=\dim \cU_t \geq \dim E_t \geq \dim E - \dim \cE \geq d+1$.
Hence $\dim \cU_t = \dim E_t$. Since $\cU_t$ is irreducible, 
this implies that 
$\cU_t=E_t$ and thus $\cU_t \subset E$.

Set $S=q(E')\subset X$. This is a set of codimension at least $2$ in $X$.
The restriction $q|_{\cU\setminus E'}:\cU\setminus E' \to X \setminus S$
is an isomorphism. The proper morphism 
$p|_{\cU\setminus E'}:\cU\setminus E'\to T\setminus p(E')$
induces a proper morphism $X\setminus S\to T\setminus p(E')$
extending $\varphi^0$.
We replace $X^0$ with  $X\setminus S$ and $T^0$ with  
$T\setminus p(E')$, obtaining a proper morphism  
$\varphi^0:X^0\to T^0$ with $\codim(X\setminus X^0)\geq 2$,
and whose general fiber is isomorphic to $\p^{d+1}$.
By shrinking $T^0$ we may assume that it is smooth.
(For this we need to remove from $T^0$ a subset of codimension 
at least $2$, and so we still have $\codim(X\setminus X^0)\geq 2$.)

Let $C$ be a curve in $T$ obtained as the intersection
of $n-d-2$ general very ample divisors.
Set $C^0=C\cap T^0$ and $X_{C^0}=(\varphi^0)^{-1}(C^0)$.
By Bertini Theorem, both $C^0$ and  $X_{C^0}$ are smooth.
Moreover, the general fiber of the induced fibration
$\varphi_{C^0}:X_{C^0}\to C^0$ is isomorphic to $\p^{d+1}$.
Since $\dim C^0 = 1$, there exists a $\varphi_{C^0}$-ample
line bundle $L$ on $X_{C^0}$ such that 
the restriction of $L$ to a general fiber of $\varphi_{C^0}$
is isomorphic to $\o_{\p^{d+1}}(1)$.
Thus we can apply \cite[Corollary 5.4]{fujita75} and conclude that 
$\varphi_{C^0}:X_{C^0}\to C^0$ is in fact a $\p^{d+1}$-bundle.
By Bertini, after removing from $T^0$ a subset of codimension at least 
$2$, we may assume that $\varphi^0:X^0\to T^0$ is in fact a 
$\p^{d+1}$-bundle.
\end{proof}


\section{Proof of Theorem~\ref{p^n}} \label{section:proof}

Let $X$ be a smooth complex projective $n$-dimensional variety.
In this section we assume that the tangent bundle $T_X$ contains 
a rank $r$ ample locally free subsheaf $E$ and prove that $X\cong \p^n$.

We begin by noticing that $X$ is uniruled.
This follows from a theorem by Miyaoka (see \cite{miyaoka} or
Shepherd-Barron's article in \cite{kollar_et_al}).
We fix a minimal covering family $H$ of rational curves on $X$ 
and set $d=\deg(f^*T_X)-2$, where $[f]$ is any member of $H$.
For a general point $x\in X$,
consider the tangent map $\tau_x:H_x\to \cC_x\subset \p(T_xX)$,
defined in section~\ref{section:C_x}.
Denote by $H_x^i$, $1\leq i\leq k$, the irreducible components of $H_x$,
and by $\pi_x^i:U_x^i\to H_x^i$ and $\eta_x^i:U_x^i\to X$
the corresponding universal family morphisms.
Denote by $\locus(H_x^i)$ the image of $\eta_x^i$, 
and by $\cC_x^i$ the image of $\tau_x|_{H_x^i}$.

We use the description of $T_{\tau_x([f])}\cC_x$ 
given in Proposition~\ref{TC_x} to
study the projective embedding
$\cC_x\subset \p(T_xX)$ for general $x\in X$.

\begin{prop} \label{C_x_is_a_linear_subspace}
Let the notation and assumptions be as above.
Then, for a general point $x\in X$, the following holds.
\begin{enumerate}
\item For every $i\in \{1,\dots,k\}$,
$\cC_x^i$  is a $d$-dimensional linear subspace of $\p(T_xX)$.
Moreover, $\p(E_x)\subset \bigcap\limits_{i=1}^k\cC_x^i$. 
\item For every $i\in \{1,\dots,k\}$,
the restriction $\tau_x|_{H_x^i}:H_x^i \to \cC_x^i$ is an isomorphism. 
As a consequence, 
all curves parametrized by $H_x$ are smooth at $x$ and 
$\eta_x^i:U_x^i\to X$ is birational onto its image for every $i\in \{1,\dots,k\}$.
\end{enumerate}
\end{prop}

\begin{proof}
Fix an irreducible component $H_x^i$ of $H_x$, and
let $[f]\in H_x^i$ be a general element. 
There is an injection of sheaves $f^*E\into f^*T_X \cong
\o_{\p^1}(2)\oplus \o_{\p^1}(1)^{\oplus d} \oplus  \o_{\p^1}^{\oplus n-d-1}$.
Since $E$ is ample, 
$f^*E$ is a subsheaf of the positive part of $f^*T_X$,
$$
f^*T_X^+ = \im [H^0(\p^1, f^*T_X(-1))\otimes \o \to f^*T_X(-1)]\otimes \o(1) 
\cong \o_{\p^1}(2)\oplus \o_{\p^1}(1)^{\oplus d},
$$
and so $E_x\subset T_xX^+_f$. 
Hence for a general element $[f]\in H_x^i$ we have
$\p(E_x)\subset \overline {T_{\tau_x([f])}\cC_x^i} \subset \p(T_xX)$.

Now we apply Lemma~\ref{cones} and conclude that each irreducible component 
$\cC_x^i$ of $\cC_x$ is a cone in $\p(T_xX)$ whose vertex contains $\p(E_x)$.

But we also know that $H_x$ is smooth and that $\tau_x:H_x\to \cC_x$ 
is the normalization.
Therefore Lemma~\ref{normalization} implies that each $\cC_x^i$ is a linear 
subspace of $\p(T_xX)$, and $\tau_x|_{H_x^i}$ is an isomorphism. 
The second part of (2) follows
from Corollary~\ref{Cx_smooth}.
\end{proof}

\begin{lemma} \label{cones}
Let $Z$ be an irreducible closed subvariety of $\p^m$. 
Assume there is a dense open subset $U$ of the smooth locus of $Z$ 
and a point $z_0\in \p^m$ such that
$z_0\in \bigcap\limits_{z\in U}T_zZ$. 
Then $Z$ is a cone in $\p^m$ and $z_0$ lies in the vertex of this cone.
\end{lemma}

\begin{proof}
%
%

We may assume that $\dim Z>0$.
Consider the projection from $z_0$, $\pi_{z_0}:Z\map \p^{m-1}$.
Since  $z_0\in T_zZ$ for general $z\in Z$, 
the tangent map to $\pi_{z_0}$ has rank $\dim Z-1$ at a general point. 
So $\pi_{z_0}$ has $1$-dimensional fibers, and thus $Z$ is a cone whose vertex
contains $z_0$. 
(Notice that in this proof we use the characteristic $0$ assumption). 
\end{proof}

\begin{lemma} \label{normalization}
If $Z$ is an irreducible cone in $\p^m$ and the normalization of $Z$ is smooth,
then $Z$ is a linear subspace of $\p^m$.
\end{lemma}

\begin{proof}
Let $x_0,\dots, x_m$ be the projective coordinates of $\p^m$.
We may assume that $Z$ is a cone with vertex $P=(0:\dots:0:1)$ over a closed
irreducible subvariety $V$ contained in the hyperplane section $(x_m=0)$ 
of $\p^m$.

Let $I_V\subset \c[x_0,\dots,x_{m-1}]$ be the homogeneous ideal defining 
$V$ in $(x_m=0)\cong \p^{m-1}$. 
By changing $m$ if necessary we may assume that $V$ 
is nondegenerate in $\p^{m-1}$.
Then $Z\setminus (x_m=0)$ has affine 
coordinate ring $S(V)=\c[x_0,\dots,x_{m-1}]/I_V$, 
and the integral closure 
of $S(V)$ is $S'=\bigoplus\limits_{l\geq 0}H^0(V,O_V(l))$.
Moreover, $S'$ can be written as $S'=\c[y_0,\dots,y_{M-1}]/I'$
for some $M>1$ and some homogeneous ideal $I'\subset \c[y_0,\dots, y_{M-1}]$.
By changing $M$ we may assume that 
$V'=\proj S'$ is nondegenerate in $\p^{M-1}$. Let 
$C(V')=\spec S' \subset \a^M$ be
the affine cone over $V'$. Then $C(V')\to Z\setminus (x_m=0)$
is the normalization morphism.

Assume $Z$ is nonlinear. Since $V$
is nondegenerate in $\p^{m-1}$, this is the same as assuming that
$I_V$ is generated by elements of degree $\geq 2$.

Now consider the inclusion of graded rings
$
S(V)=\c[x_0,\dots,x_{m-1}]/I_V\into S'=\c[y_0,\dots,y_{M-1}]/I',
$ 
and denote by $\varphi_i$ the image of $x_i$ in $S'_1$.
Since $x_0,\dots, x_{m-1}$ are linearly independent in $S(V)_1$,
$\varphi_0,\dots, \varphi_{m-1}$ are linearly independent in $S'_1$.
But then $M \geq m > \dim Z = \dim C(V')$.
Since we assume that $V'$ is nondegenerate,
this implies that $C(V')$ is a nonlinear cone,
and hence not smooth.
\end{proof}

The next step is to prove that $H_x$ is in fact irreducible for general $x\in X$.
The idea is to produce a curve $C$ through $x$ such that, for every $i\in \{1,\dots,k\}$,
there exists an element $[f_i]\in H_x^i$ parametrizing $C$.
Since $C$ is smooth at $x$ (by Proposition~\ref{C_x_is_a_linear_subspace}(2))
there exists a unique 
point in $H_x$ parametrizing $C$, and 
$H_x$ must be irreducible.

By Proposition~\ref{C_x_is_a_linear_subspace}(1), $\cC_x$ is 
the union of linear subspaces of $\p(T_xX)$ for general $x\in X$.
Fix $x\in X$ and let $i\in \{1,\dots,k\}$. Set $Y_i=\locus(H_{x}^i)$.
By Lemma~\ref{leaf_p^n}, 
the normalization of $Y_i$ is isomorphic to $\p^{d+1}$.
Under this isomorphism, the rational curves on $Y_i$ parametrized by $H_{x}^i$
come from the lines on $\p^{d+1}$ passing through a fixed point $\tilde x_i\in \p^{d+1}$.
Let $n_i:\p^{d+1}\to Y_i$ be the normalization morphism.

We claim that the vector bundle $E|_{Y_i}$ pulls back to a subsheaf of $T_{\p^{d+1}}$.
Indeed,
the injection $E\into T_X$ induces a map $\Omega_X \to E^{\vee}$
of maximal rank.
By restricting to $Y_i$, we get a map $\Omega_X|_{Y_i}\to E^{\vee}|_{Y_i}$
of maximal rank.
This map factors through  $\Omega_{Y_i}\to E^{\vee}|_{Y_i}$. (This is because 
the composite map $I_{Y_i}/I^2_{Y_i}\to \Omega_X|_{Y_i}\to E^{\vee}|_{Y_i}$
vanishes identically.)
Lemma~\ref{L_lifts} below asserts that there is a map 
$\Omega_{\p^{d+1}}\to n_i^*E^{\vee}|_{Y_i}$ factoring 
$n_i^*\Omega_{Y_i}\to n_i^*E^{\vee}|_{Y_i}$. By dualizing we get 
a sheaf injection $n_i^*E|_{Y_i} \into T_{\p^{d+1}}$.

Thus, $n_i^*E|_{Y_i}$ is an ample vector bundle on $\p^{d+1}$
that is a subsheaf of $T_{\p^{d+1}}$.  
So either $n_i^*E|_{Y_i}\cong \o_{\p^{d+1}}(1)^{\oplus r}$,
or $r=d+1$ and $n_i^*E|_{Y_i}\cong T_{\p^{d+1}}$.

In any case, there exists an $r$-dimensional linear subspace $M_i$ of $\p^{d+1}$,
passing through $\tilde x_i$, for which 
$(n_i^*E|_{Y_i})|_{M_i}\into T_{\p^{d+1}}|_{M_i}$ factors through 
$T_{M_i}\into T_{\p^{d+1}}|_{M_i}$.
Set $Z=n_i(M_i)\subset Y_i\subset X$.
Then $Z$ is an $r$-dimensional subvariety of $X$ containing $x$ 
and tangent to $E$ along its smooth locus. 
By Frobenius' Theorem, $Z$ is the unique subvariety of $X$ with these properties. 
Hence $Z=n_i(M_i)$ for every $i\in \{1,\dots,k\}$.
Let $C$ be the image in $Z$ of a line through $\tilde x_i$ on 
$M_i$ for some $i$. 
Then, for every $i$, $C$ comes from a line through $\tilde x_i$ on 
$M_i$, and thus 
there exists an element $[f_i]\in H_x^i$ parametrizing $C$.
This shows that $H_x$ is irreducible as we noted above.

\begin{lemma} \label{L_lifts}
Let $Z$ be a variety, $n:\tilde Z \to Z$ the normalization, and 
$E$ a vector bundle on $Z$. 
Assume there is a map $\Omega_Z\to E^{\vee}$. 
Then the induced map $n^*\Omega_Z\to n^*E^{\vee}$ factors through
$n^*\Omega_Z \to \Omega_{\tilde Z}$.
\end{lemma}

\begin{proof}
This result follows from a theorem by Seidenberg, which asserts that a 
derivation of an integral domain over a ground field of characteristic $0$
extends to its normalization.

Let $U=\spec A\subset Z$ be an affine open subset over which $E$ is trivial,
and fix an isomorphism $E|_U\cong \o_U^{\oplus r}$, where $r=\rank E$.

Let $\tilde A$ be the integral closure of $A$.

The restricted map $\Omega_U\to E|_U\cong \o_U^{\oplus r}$ 
induces a homomorphism
$\Omega_A\to A^{\oplus r}$.
By composing with the $r$ natural projections,
$p_i:A^{\oplus r}\to A$, $1\leq i \leq r$, we obtain $r$ 
derivations $D_i:A\to A$, $1\leq i \leq r$. 
By the main theorem in \cite{seidenberg}, 
each of these derivations extends uniquely
to a derivation $\tilde D_i:\tilde A\to \tilde A$. 
Such $\tilde D_i$'s determine a homomorphism 
$\Omega_{\tilde A}\to \tilde A^{\oplus r}$
extending $\Omega_A\to A^{\oplus r}$.   
\end{proof}

Theorem~\ref{Cx_fibration} yields 
a dense open subset $X^0$ of $X$ and a $\p^{d+1}$-bundle 
$\varphi^0:X^0\to T^0$.
Since an ample vector bundle of rank $r$ on a rational curve has degree
at least $r$,
either $H$ is an unsplit family, or
$d=0$, $r=1$, and $f^*E\cong \o_{\p^1}(2)$
for every $[f]\in H$.
We analyse these two cases separately.

\medskip 
\noindent
{\bf Case 1 ($H$ is an unsplit family). }

In this case we can apply  Theorem~\ref{extending_in_codim_1}
and assume that $T^0$ is smooth and $\codim(X\setminus X^0)\geq 2$.

Suppose that $\dim T^0 > 0$. Let $C'\subset X^0$ 
be a general smooth projective curve
such that $C=\varphi^0(C')$ is also a smooth projective curve.
(Such a curve exists because $X$ and $T^0$ are smooth and 
$X\setminus X^0$ has codimension at least $2$.)
Then $X_C:=(\varphi^0)^{-1}C\to C$ is a $\p^{d+1}$-bundle.
Since $C'$ is general, there is a sheaf inclusion 
$E|_{X_C}\into T_X|_{X_C}$.
For general $x\in X_C$ we have 
$E_x\subset (T_{X_C/C})_x\subset T_xX$.
The cokernel of the map $T_{X_C/C}\into T_X|_{X_C}$
is torsion free.
Hence $E|_{X_C}$ is in fact a subsheaf of the relative 
tangent sheaf  $T_{X_C/C}$. But this contradicts 
Lemma~\ref{campana-peternell} below, due to Campana and Peternell. 
Therefore $T^0$ is a point, $X\cong \p^n$ and under this 
isomorphism either $E=T_{\p^n}$ or $E\cong \o_{\p^n}(1)^{\oplus r}$.

\begin{lemma}[{\cite[Lemma 1.2]{campana_peternell}}]
\label{campana-peternell}
Let $T$ be a smooth complex projective variety of positive dimension, 
$E$ a vector bundle of rank $k+1$ on $T$, and 
$X=\PP(E)\to T$ the corresponding $\p^k$-bundle.
Then the relative tangent sheaf $T_{X/T}$ does not contain 
any ample locally free subsheaf.
\end{lemma}

\noindent
{\bf Case 2 ($H$ is not proper, $d=0$, $r=1$ and $f^*E\cong \o_{\p^1}(2)$
                              for every $[f]\in H$).}

We have $\dim H>0$. Let $C'\subset H$ be a general curve. 
Let $C$ be the normalization of the closure of $C'$ 
in $\chow(X)$. 
(Notice that, since $H$ is not proper, some points of $C$ may
parametrize nonintegral curves.)
Let $S$ be the normalization of the universal family 
over $C$ and denote by 
$p:S\to C$ and $n:S\to X$ the universal family morphisms.

Let $S'=n(S)\subset X$. Then $n:S\to S'$ is birational.
Since $\dim H_x=0$ for general $x\in X$, $n$ does not contract
any curve dominating $C$. Neither does it contract any curve
contained in a fiber of $p$. Hence $n$ is the normalization.
By Lemma~\ref{L_lifts}, the injection $E|_{S'}\into \Omega_{S'}^{\vee}$
lifts to an injection $n^*E\into \Omega_S^{\vee}$.
For convenience set $L=n^*E$.

The idea is to reach a contradiction as follows.
We look at the minimal resolution of $S$
and contract the $(-1)$-curves that do not dominate $C$.
In this way we obtain a $\p^1$-bundle over $C$.
We show that $L$ induces an ample line bundle on the resulting
$\p^1$-bundle that is a subsheaf of the relative tangent sheaf.
But this is impossible by Lemma~\ref{campana-peternell}.

So let $r:Y\to S$ be the minimal resolution of $S$ and set $L_Y=r^*L$ 
(notice that $L_Y$ is an ample line bundle on $Y$).
By \cite[Proposition 1.2]{burns_wahl}, there is a natural isomorphism 
$r_*T_Y\xrightarrow{\cong} \Omega_S^{\vee}$.
Therefore, from the natural isomorphism
$\Hom_Y(L_Y,T_Y)\cong \Hom_S(L,r_*T_Y)$
(see \cite[II.5]{hartshorne}), we see that 
the map $L\to \Omega_S^{\vee}$
lifts to an injection $L_Y\into T_Y$.

The induced morphism $p_Y:Y\to C$ can be obtained 
from a suitable $\p^1$-bundle $p_Z:Z\to C$
by a composition of blowups, $q:Y\to Z$.
Set $L_Z=q_*L_Y$.

By pushing forward to $Z$ and applying the projection formula, we see that
the inclusion $L_Y\into T_Y$ induces an inclusion $L_Z\into T_Z$.

To show that $L_Z$ is an ample line bundle, it is enough to
assume that $q:Y\to Z$ is the inverse of a single blowup 
(then use induction on the number of blowups). 
First note that $L_Z$ is in fact a line bundle on $Y$ 
(it is reflexive except possibly at finitely many points, and
hence reflexive). 
So we can write $L_Y=q^*L_Z+aD$, where $D$ is the exceptional 
curve, and $a=-aD^2=-L_Y\cdot D$ is a negative integer.
The ampleness of $L_Z$ then follows from Nakai's criterion.

For any fiber $F$ of $p_Z$ we have $L_Z\cdot F>0$. Hence, for a general 
fiber $F\cong \p^1$, the map $L_Z|_F\into T_Z|_F\cong \o_{\p^1}(2)\oplus \o_{\p^1}$ 
factors through 
$\o_{\p^1}(2)\cong T_F\into T_Z|_F$. 
Since the cokernel of the map $T_{Z/C}\into T_Z$ is 
torsion free, this implies that 
there is an inclusion $L_Z\into T_{Z/C}$ factoring 
$L_Z\into T_Z$.

We have shown that $L_Z$ is an ample line bundle on $Z$
that injects into $T_{Z/C}$.
This is a contradiction as we noted above.
So case 2 does not occur.

\medskip
\noindent {\it Acknowledgements. }
Most of the content of this article appears in my Ph.D. thesis.
I am deeply grateful to my advisor, J\'anos Koll\'ar, for his 
guidance, attention, many ideas and suggestions.
Financial support during graduate school was provided by 
CNPq (Conselho Nacional de Desenvolvimento Cient\'\i fico e Tecnol\'ogico --- Brazil)
and Princeton University.
This research was partially conducted during the period I was employed by the 
Clay Mathematics Institute as a Liftoff Fellow.

\bibliographystyle{amsalpha}
\bibliography{carolina_paper}

\providecommand{\bysame}{\leavevmode\hbox to3em{\hrulefill}\thinspace}
\providecommand{\MR}{\relax\ifhmode\unskip\space\fi MR }
\providecommand{\MRhref}[2]{%
  \href{http://www.ams.org/mathscinet-getitem?mr=#1}{#2}
}
\providecommand{\href}[2]{#2}
\begin{thebibliography}{{Kol}92}

\bibitem[AW01]{andreatta_wisniewski}
M.~Andreatta and J.~A. Wi{\'s}niewski, \emph{On manifolds whose tangent bundle
  contains an ample subbundle}, Invent. Math. \textbf{146} (2001), no.~1,
  209--217.

\bibitem[BW74]{burns_wahl}
D.~M. Burns, Jr. and J.~M. Wahl, \emph{Local contributions to global
  deformations of surfaces}, Invent. Math. \textbf{26} (1974), 67--88.

\bibitem[CP98]{campana_peternell}
F.~Campana and T.~Peternell, \emph{Rational curves and ampleness properties of
  the tangent bundle of algebraic varieties}, Manuscripta Math. \textbf{97}
  (1998), no.~1, 59--74.

\bibitem[Dru04]{druel}
S.~Druel, \emph{Caract\'erisation de l'espace projectif}, Manuscripta Math.
  \textbf{115} (2004), no.~1, 19--30.

\bibitem[Fuj75]{fujita75}
T.~Fujita, \emph{On the structure of polarized varieties with {$\Delta
  $}-genera zero}, J. Fac. Sci. Univ. Tokyo Sect. IA Math. \textbf{22} (1975),
  103--115.

\bibitem[Har77]{hartshorne}
R.~Hartshorne, \emph{Algebraic geometry}, Springer-Verlag, New York, 1977,
  Graduate Texts in Mathematics, No. 52.

\bibitem[HM04]{hwang_mok_birationality}
J.-M. Hwang and N.~Mok, \emph{Birationality of the tangent map for minimal
  rational curves}, Asian J. Math. \textbf{8} (2004), no.~1, 51--64.

\bibitem[Hwa01]{hwang}
J.-M. Hwang, \emph{Geometry of minimal rational curves on {F}ano manifolds},
  School on Vanishing Theorems and Effective Results in Algebraic Geometry
  (Trieste, 2000), ICTP Lect. Notes, vol.~6, Abdus Salam Int. Cent. Theoret.
  Phys., Trieste, 2001, pp.~335--393.

\bibitem[Keb02]{kebekus}
S.~Kebekus, \emph{Families of singular rational curves}, J. Algebraic Geom.
  \textbf{11} (2002), no.~2, 245--256.

\bibitem[{Kol}92]{kollar_et_al}
{Koll{\'a}r, J. et al}, \emph{{F}lips and abundance for algebraic threefolds},
  Soci\'et\'e Math\'ematique de France, Paris, 1992, Ast\'erisque No. 211
  (1992).

\bibitem[Kol96]{kollar}
J.~Koll{\'a}r, \emph{Rational curves on algebraic varieties}, Ergebnisse der
  Mathematik und ihrer Grenzgebiete, vol.~32, Springer-Verlag, Berlin, 1996.

\bibitem[Miy87]{miyaoka}
Y.~Miyaoka, \emph{The {C}hern classes and {K}odaira dimension of a minimal
  variety}, Algebraic geometry, Sendai, 1985, Adv. Stud. Pure Math., vol.~10,
  North-Holland, Amsterdam, 1987, pp.~449--476.

\bibitem[Mor79]{mori79}
S.~Mori, \emph{Projective manifolds with ample tangent bundles}, Ann. of Math.
  (2) \textbf{110} (1979), no.~3, 593--606.

\bibitem[Sei66]{seidenberg}
A.~Seidenberg, \emph{Derivations and integral closure}, Pacific J. Math.
  \textbf{16} (1966), 167--173.

\bibitem[Wah83]{wahl}
J.~M. Wahl, \emph{A cohomological characterization of {${\bf P}\sp{n}$}},
  Invent. Math. \textbf{72} (1983), no.~2, 315--322.

\end{thebibliography}

\end{document}